\documentclass[11pt]{article}
\usepackage{graphicx,psfrag}
\usepackage{amssymb,amsfonts,amsthm}
\usepackage{amsmath,amsthm,amssymb}
\usepackage{amsfonts}

\newcommand{\la}{\langle}
\newcommand{\ra}{\rangle}

\newtheorem{theorem}{Theorem}
\newtheorem{lemma}[theorem]{Lemma}

\newtheorem{cor}[theorem]{Corollary}
\newtheorem{example}[theorem]{Example}
\newtheorem{proposition}[theorem]{Proposition}
\theoremstyle{definition}

\newtheorem{prob}[theorem]{Problem}

\newcommand{\trace}{{\mathrm{trace}}}
\newcommand{\HNN}{\mathrm{HNN}}

\newcommand {\C}{\mathbb{C}} 

\newcommand {\iv}{^{-1}}

\begin{document}

\title{Non-linear residually finite groups}
\author{Cornelia Dru\c{t}u and Mark Sapir\thanks{The first author
is grateful to the CNRS of France for granting her the
d\'el\'egation CNRS status in the 2003-2004 academic year. The
research of the second author was supported in part by the NSF
grant DMS 0072307.}}
\date{}
\maketitle

\begin{abstract} We give the first example
of a non-linear residually finite $1$-related group: $\la a, t\mid
a^{t^2}=a^2\ra$.
\end{abstract}

\section{Non-linear groups}

Let $\phi$ be an injective endomorphism of a group $G$. Then the
HNN extension $$\mathrm{HNN}_\phi(G)=\la G,t\mid tgt\iv=\phi(g),
g\in G\ra$$ is called an {\em ascending HNN extension} of $G$ (or
the {\em mapping torus of the endomorphism $\phi$}). In
particular, the ascending HNN extensions of free groups of finite
rank are simply the groups given by presentations $\la
x_1,...,x_n, t\mid tx_it\iv=w_i, i=1,...,n\ra$, where $w_1,...,w_n$
are words in $x_1,...,x_n$ generating a free subgroup of rank $n$.

In \cite{BS}, Borisov and Sapir proved that all ascending HNN
extensions of linear groups are residually finite. After \cite{BS},
the question of linearity of these groups became very interesting.
This question is especially interesting for ascending HNN extensions
of free groups because most of these groups are hyperbolic
\cite{Kap}, and because calculations show that at least $99.6\%$ of
all 1-related groups are ascending HNN extensions of free groups
\cite{BS}.

 Let
$H=\mathrm{HNN}_\phi(F_n)$ be an ascending HNN extension of a free
group. If $n=1$ then $H$ is a Baumslag-Solitar group $BS(m,1)$, so
it is inside $\mathrm{SL}_2(\mathbb{Q})$. If $n=2$ and $\phi$ is an
automorphism then the linearity of $H$ follows from the linearity of
$\mathrm{Aut}(F_2)$. The linearity of $\mathrm{Aut}(F_2)$ follows
from two facts: Dyer, Formanek and Grossman \cite{DFG} reduced the
linearity of $\mathrm{Aut}(F_2)$ to the linearity of the braid group
$B_4$;  the linearity of $B_4$ was proved by Krammer \cite{Kra}.

It is known that in the case when $\phi$ is not an automorphism,
the situation is different.

\begin{proposition}[Wehrfritz, \cite{W}, Corollary 2.4]\label{main}
The group $\la a,b, t\mid tat\iv =a^k, tbt\iv=b^l\ra$, with
$k,l\not\in\{1,-1\}$, is not linear.
\end{proposition}

The proof\footnote{In the first version of this paper that appeared
in the arXiv, we gave a complete proof of this statement because we
were unaware of \cite{W}. We are grateful to Professor Raptis for
providing this reference.} in \cite{W} uses the action of
$\mathrm{SL}_n(K)$ on the Lie algebra $\mathrm{sl}_n(K)$ to deduce
that if $a,b,t$ are matrices satisfying the relations of the group then
some powers of $a$ and $b$ generate a nilpotent subgroup. One can
also prove this statement by using the fact that if $a, t$ are
matrices with complex entries such that $tat\iv=a^k$, $|k|>1$, then
$a^{k-1}\in U(t\iv )=\{x\mid \lim_{n\to\infty}t^{-n}xt^n = {\bf
1}\}$ and the well known fact that $U(t\iv )$ is a
nilpotent group for every matrix $t$.

The following lemma is useful when dealing with ascending HNN
extensions of groups.

\begin{lemma}\label{lm1}
Let $\phi$ be an injective endomorphism of a group $G$.
Suppose that $\phi^k$ is not an inner automorphism of $G$
for any $k\ne 0$. Then a homomorphism $\gamma$ of $H=\HNN_\phi(G)$
is injective if and only if the restriction of $\gamma$ on $G$ is
injective.
\end{lemma}

\proof This follows immediately from the fact \cite{BS} that every
element of $H$ is of the form $t^{-p}wt^q$, where $p,q\ge0$, $w\in
G$. Indeed suppose that $\gamma$ is injective on $G$ but
$\gamma(t^{-p}wt^q)=1$ where $t^{-p}wt^q\ne 1$. Then
$\gamma(wt^{q-p})=1$, and so $q\ne p$. Hence
$\gamma(w)=\gamma(t^{p-q})$. Let $m=p-q$. We can assume that $m>0$,
otherwise replace $w$ by $w\iv$. This implies that for every $u\in
G$, $\gamma(t^mut^{-m})=\gamma(wuw\iv)$. The injectivity of the restriction of
$\gamma$ on $G$ then implies that $t^mut^{-m}=wuw\iv$. Since $t^m
ut^{-m}=\phi^m(u)$, we get $\phi^m(u)=wuw\iv$. Hence $\phi^m$
is the inner automorphism induced by $w$, a contradiction.
\endproof

\begin{cor}\label{one}
The group $H=\la a,b, t\mid tat\iv =a^k, tbt\iv=b^l\ra$ is linear if
and only if $k, l\in\{1,-1\}$.
\end{cor}

\proof If $k,l\in \{1,-1\}$ then $\phi$ is an automorphism and $H$
is linear, say, by the results of \cite{DFG} and \cite{Kra} cited
above (one can also use the fact that the group has a subgroup of
finite index isomorphic to $F_2\times\mathbb{Z}$). If both
$k,l\not\in\{1,-1\}$, we can apply Proposition \ref{main}. It
remains to consider the case $k\not\in\{1,-1\}$, $l\in\{1,-1\}$.
Then it is easy to see by Lemma \ref{lm1} that the subgroup $\la a,
bab\iv, t^2\ra$ is isomorphic to the group $\la x, y, t\mid
txt\iv=x^{k^2}, tyt\iv=y^{k^2}\ra$, and so it is not linear by
Proposition \ref{main}. Hence the group $H$ is not linear as well.
\endproof

Not much is known about the linearity of 1-related groups. Note only
that all residually finite Baumslag-Solitar groups (i.e. HNN
extensions of cyclic groups) \cite{Me} are linear \cite{Vo}.

The following theorem gives the first example of a non-linear
residually finite 1-related group.

\begin{theorem}\label{main1} The group $\la a, t\mid
t^2at^{-2}=a^2\ra$ is residually finite but not linear.
\end{theorem}

\proof Using Magnus rewriting procedure, this group can be
represented as an HNN extension $\la a,b,t\mid tat\iv = b,
tbt\iv=a^2\ra$, so it is residually finite by \cite{BS}. The
subgroup of that group generated by $\{a,b,t^2\}$ is isomorphic (by
Lemma \ref{lm1}) to $\la a,b, t\mid tat\iv =a^2, tbt\iv=b^2\ra$
which is not linear by Proposition \ref{main}. The isomorphism takes
$a$ to $a$, $b$ to $b$, $t$ to $t^2$.
\endproof

\begin{prob}\label{p1} Is it true that $\mathrm{HNN}_\phi(F_n)$ is always linear if
$\phi$ is an automorphism?
\end{prob}

\begin{prob}\label{p2} Are there hyperbolic non-linear ascending HNN
extensions of free groups? In particular, is the group $\la a,b,
t\mid tat\iv =ab, tbt\iv=ba\ra$ linear (the fact that this group
is hyperbolic follows from \cite{Kap})?
\end{prob}


The group $\la a,b,t\mid tat\iv=ab, tbt\iv=ba\ra$ is actually a
1-related group $\la a,t\mid [[a,t],t]=a\ra$. The fact that this
group does not have a faithful 2-dimensional representation follows
from \cite{FLR}. Moreover, results of \cite{FLR} (and prior results
of Magnus \cite{M}) imply that most 1-related groups do not have
faithful 2-dimensional representations.

We conjecture that the answer to Problem \ref{p2} is that most
groups $HNN_\phi(F_k)$, $k\ge 2$, are non-linear provided $\phi$ is
not an automorphism. In particular there are many non-linear
hyperbolic groups among these HNN extensions.

Note that there are examples of non-linear hyperbolic groups yet,
although M.Kapovich \cite{Kap1} has an example of a hyperbolic group
which does not have faithful real linear
representations.\footnote{Added 10/12/04: It is easy to show using
Kapovich's argument that the group does not have faithful
representations over any field.}

Proposition \ref{main} and Corollary \ref{one} give examples of
non-linear ascending HNN extensions of linear groups. Non-ascending
HNN extensions with this property are much easier to find: some of
these HNN extensions are not even residually finite (say, the
Baumslag-Solitar groups $\la a,t \mid ta^2t\iv=a^3\ra$). Residually
finite non-linear HNN extensions of linear groups were constructed
in particular by Formanek and Procesi \cite{FP}. They proved that
the HNN extension of the direct product $F_k \times F_k$, where one
of the associated subgroups is the diagonal and the other one is one
of the factors, is residually finite but not linear. This was the
main ingredient in the proof in \cite{FP} of the non-linearity of
$\mathrm{Aut}(F_n)$, $n\ge 3$.

\section{Representations in $\mathrm{SL}_2(\mathbb{C})$}

 By Lemma \ref{lm1}, finding a copy of $H=\HNN_\phi(F_k)$ in ${\mathrm
SL}_n(K)$ amounts to finding a $k$-tuple of matrices $(A_1,...,A_k)$
that freely generate a free subgroup, and which is a conjugate of
the $k$-tuple $(\phi(A_1),...,\phi(A_k))$. In the case when $k=2$,
$n=2$, $K=\mathbb{C}$ one can use the fact that conjugacy of two
pairs of $2$ by $2$ matrices $(A,B)$, $(C,D)$ implies the system of
equalities $\trace(A)=\trace(C), \trace(B)=\trace(D),
\trace(AB)=\trace(CD)$. The converse statement is ``almost" true:
one needs to exclude the case when
\begin{equation}\label{3}
\trace([A,B])=\trace(A)^2+\trace(B)^2+\trace(AB)^2-\trace(A)\trace(B)\trace(AB)-2=2
\end{equation} (in that case $A,B$ generate a solvable group
\cite{Bow}). Using the fact that for every word $u=u(A,B)$ in
matrices $A,B \in \mathrm{SL}_2(\mathbb{C})$, $\trace(u)$ can be
expressed as a polynomial in $\trace(A)$, $\trace(B)$, $\trace(AB)$,
we get a system of three equations with three unknowns. The
corresponding algebraic variety will be called the {\em trace
variety} of the group $\HNN_\phi(F_k)$.

In most cases that we considered, the trace variety was
0-dimensional. But the next example shows that the trace variety
may have dimension $\ge 1$ and the group still may not have a
faithful 2-dimensional representation.

\begin{example} The group $H=\la a, b, t| tat\iv=a, tbt\iv=[a,b]\ra$ does not have a faithful
2-dimensional representation. The trace variety of this group is a
union of two curves, but it consists of non-faithful
representations.
\end{example}

Consider any representation of $H$ in $\mathrm{SL}_2(\mathbb{C})$.
So we assume that $a,b$ are $2$ by $2$ matrices with determinant
$1$. Let us denote $\trace(a)=x, \trace(b)=y, \trace(ab)=z$. It is
easy to see using
\begin{equation}\trace(BA^2C)=\trace(A)\trace(BAC)-\trace(BC)\label{23}\end{equation}
(this is essentially the Cayley-Hamilton theorem for matrices in
$\mathrm{SL}_2$) that we have the following system of equations:
$$\left\{\begin{array}{l}
x=x\\
y= \trace(aba\iv b\iv)=-2+x^2+y^2+z^2-xyz\\
z=\trace(a[a,b])=_{\mathrm{by (\ref{23})}} x\cdot\trace(aba\iv
b\iv)-\trace(ba\iv b\iv)\\ \hskip .12 in  = x\cdot\trace(aba\iv
b\iv)-x=xy-x
\end{array}\right.$$

Plugging $z=xy-x$ into the second equation, and solving for $y$, we
get $y=2, x=z$ or $y=x^2-1,z=x^3-2x$. Thus the trace variety is a
union of two curves. If $y=2, x=z$ then $\la a,b\ra$ is solvable by
(\ref{3}), so the representation is not faithful. Now let $y=x^2-1,
z=x^3-2x$. Consider the word $w=ab\iv a\iv ba\iv b\iv a$. It is not
difficult to compute the corresponding trace polynomial:
$$\trace(w)=-3y-4xz+5yx^2+xz^3-2yx^2z^2+yz^2+y^3-y^3x^2+y^2x^3z+x^3z-yx^4.$$
If we plug in $y=x^2-1, z=x^3-2x$ into this polynomial, we get $2$.
Similarly, the trace polynomial of the word $wa$ is

$$x^4y^2z-x^5y-x^3y^3-2x^3yz^2+x^4z-x^2y^2z+x^2z^3+6x^3y+2xy^3+3xyz^2-5x^2z-y^2z-z^3-7xy+3z.$$
If we plug in $y=x^2-1, z=x^3-2x$, we get $x$. Hence $\trace(w)=2,
\trace(wa)=\trace(a)=x$, whence $\trace([w,a])=2$ (see (\ref{3}))
and $w$ and $a$ generate a solvable subgroup. Therefore for every
value of $x$ the corresponding representation of the group $H$ is
not faithful. (In fact it is not difficult to show that the relation
$(a^2b)^3=1$ also holds, so in this case $\la a,b\ra$ has torsion.)

Similarly the trace variety of the group $\la a,b,t\mid tat\iv=a,
tbt\iv=(ba)b(ba)\iv\ra$ is two-dimensional, but this group does not
have faithful representations in $\mathrm{SL}_2(\mathbb{C})$.

In fact we do not know the answer to the following question.

\begin{prob} Are there any ascending HNN extensions of $F_k$, $k>1$, which
have faithful $2$-dimensional complex representations? In
particular, are there free non-cyclic subgroups in
$\mathrm{SL}_2(\mathbb{C})$ which are conjugate inside
$\mathrm{SL}_2(\mathbb{C})$ to their proper subgroups?
\end{prob}

\begin{minipage}[t]{2.9 in}
\noindent Cornelia Dru\c tu\\ Department of Mathematics\\
University of Lille 1\\ and UMR CNRS 8524 \\ Cornelia.Drutu@math.univ-lille1.fr\\
\end{minipage}
\begin{minipage}[t]{2.6 in}
\noindent Mark V. Sapir\\ Department of Mathematics\\
Vanderbilt University\\
msapir@math.vanderbilt.edu\\
\end{minipage}

\begin{thebibliography}{ECHPT}
\label{bibbb}

\bibitem[BS]{BS} A. Borisov, M. Sapir. Polynomial maps over finite fields and residual finiteness of
mapping tori of group endomorphisms, arXiv:math.GR/0309121.

\bibitem[Bow]{Bow} B.H. Bowditch. Markoff triples and quasi-Fuchsian groups. Proc.
London Math. Soc. (3) 77 (1998), no. 3, 697--736.

\bibitem[DFG]{DFG} J. L. Dyer, E. Formanek, E. K. Grossman. On the
linearity of automorphism groups of free groups. Arch. Math.
(Basel) 38 (1982), no. 5, 404--409.


\bibitem[FLR]{FLR} B. Fine, F. Levin, G. Rosenberger.
Faithful complex representations of one relator groups. New
Zealand J. Math. 26 (1997), no. 1, 45--52.

\bibitem[FP]{FP} E. Formanek, C. Procesi. The
automorphism group of a free group is not linear. J. Algebra 149
(1992), no. 2, 494--499.



\bibitem[I.Kap]{Kap} I. Kapovich. Mapping tori of endomorphisms of free
groups. Comm. Algebra 28 (2000), no. 6, 2895--2917.

\bibitem[M.Kap]{Kap1} M. Kapovich, Representations of polygons of
finite groups, preprint,
http://www.math.ucdavis.edu/$\sim$kapovich/EPR/gr4.ps .


\bibitem[Kra]{Kra} D. Krammer. The braid group $B\sb 4$ is
  linear. Invent. Math. 142 (2000), no. 3, 451--486.

\bibitem[M]{M} W. Magnus. Two generator subgroups of ${\rm PSL}$ $(2,\C)$.
Nachr. Akad. Wiss. G\"ottingen Math.-Phys. Kl. II 1975, no. 7, 81--94.

\bibitem[Me]{Me} S. Meskin. Nonresidually finite one-relator groups. Trans.
Amer. Math. Soc. 164 (1972), 105--114.

\bibitem[Vo]{Vo} R.T. Vol'vachev. Linear representation of certain groups with
one relation. Vestsi Akad. Navuk BSSR Ser. Fi z.-Mat. Navuk 1985,
no. 6, 3--11, 124.

\bibitem[W]{W} B.A.F. Wehrfritz.
Generalized free products of linear groups. Proc. London Math. Soc.
(3) 27 (1973), 402--424.

\end{thebibliography}
\end{document}